\documentclass[20pt,twoside,openany,reqno]{amsart}
 \usepackage{bbm}
 \usepackage{mathrsfs}
\input amssymb.sty
\usepackage{amsmath}
\usepackage[all,poly,knot]{xy}
\usepackage{amsthm}
\usepackage{amsfonts}
\usepackage{latexsym}
\usepackage{amscd}
\newcommand{\pf}{\noindent\begin {proof}}
\newcommand{\epf}{\end{proof}}

\def\bc{\begin{center}}
\def\ec{\end{center}}

\begin{document}

{\centerline{\LARGE {\bf Remarks on Murre's conjecture on Chow
groups$^{*}$}}

\bigskip

\bigskip

{\centerline{by}}

\bigskip

{\centerline{ { {\large K}EJIAN {\large X}U} \ {\small AND} \ {
{\large Z}E {\large X}U}}}

\bigskip

\bigskip

\begin{figure}[b]
\rule[-2.5truemm]{5cm}{0.1truemm}\\[2mm]
{\small
\begin{tabular}{ll}
$^{\ast}$ This research is supported by National Natural Foundation
of China (10871106).
\end{tabular}}
\end{figure}

{\centerline{\bf Abstract.}}

\bigskip
For certain product varieties, Murre's conjecture on Chow groups is
investigated. \ \ \ \

In particular, it is proved that Murre's conjecture (B) is true for
two kinds of four-\ \ \ \ \ \ \ \

folds. Precisely, if $C$ is a curve and $X$ is an elliptic modular
threefold over $k$ (an \ \ \ \

algebraically closed field of characteristic $0$) or  an abelian
variety of dimension 3, \ \

then Murre's conjecture (B) is true for the fourfold $X\times C.$}

\bigskip

{{\it Key Words:} \  motivic decomposition, Chow group, curve,
abelian variety, elliptic

modular threefold}

\bigskip
{{\it Mathematics Subject Classification 2010:} 14C25

\bigskip

\bigskip

\bigskip

{\centerline{\bf 1. Introduction}}

\bigskip

\noindent We will work with the category $\mathscr{V}_{k}$ of smooth projective varieties over a field $k$. Let $X\in \mathscr{V}_{k}$ be irreducible and of dimension
$d$. Let $H(X):=H_{et}^{*}(\overline{X}, \mathbb{Q}_{l})$ be the
$l$-adic cohomology groups over a (fixed) algebraic closure $\overline{k}$
of $k$, where $\overline{X}=X\times_{k} \text{Spec}(\overline{k})$
and $l\neq \text{ch}(k)$ is a prime, and let  $\text{cl}_{X}:
Z^{i}(X)\rightarrow H^{2i}(X)$ be the  cycle map associated to
$H(X),$ where $Z^{i}(X)$ is the group of algebraic cycles of
codimension $i$ of $X.$  We have the well-known K\"{u}nneth formula:
$$H^{2d}(X\times X)\simeq \bigoplus_{i=0}^{2d} H^{2d-i}(X) \otimes H^{i}
(X).$$
 Let $\Delta_{X} \subseteq X\times X$ be the diagonal. Then
$\text{cl}_{X\times X}(\Delta_{X})$ has the K\"{u}nneth
decomposition:
$$\text{cl}_{X\times X}(\Delta_{X})=\pi_{0}^{\text{hom}}+\pi_{1}^{\text{hom}}+\ldots+\pi_{2d}^{\text{hom}},$$
where $\pi _{i}^{\text{hom}}\in H^{2d-i}(X) \otimes H^{i} (X)$
is the $i$-th K\"{u}nneth component.

Let $A_{\text{num}}^{j}(X)$ (resp.
$A_{\text{rat}}^{j}(X)=\text{CH}^{j}(X)$) be the groups of algebraic
cycles of codimension $j$ modulo the numerical equivalence (resp.
rational equivalence). Grothendieck's Lefschetz standard conjecture implies the  $\pi _{i}^{\text{hom}}$ are all algebraic (i.e., they are all in the image of the cycle map). Assuming additionally the conjecture that the homological equivalence coincides with the numerical equivalence ([13]), the diagonal (modulo the numerical equivalence) has a canonical decomposition into a sum of orthogonal idempotents (also called projectors)
$$\Delta_{X}=\pi_{0}^{\text{num}}+\pi_{1}^{\text{num}}+\ldots+\pi_{2d}^{\text{num}}.\eqno (1)$$ in the correspondence ring $A_{\text{num}}^{d}(X\times
X)\otimes_{\mathbb{Z}}\mathbb{Q}$. Then, in the category of Grothendieck
motives $\mathscr{M}_{k}^{\text{num}}$ ([13]) (w.r.t. the numerical
equivalence), the motive $h(X)\in \mathscr{M}_{k}^{\text{num}}$ has a
canonical decomposition
$$h(X)=h^{0}(X)\oplus h^{1}(X)\oplus \ldots \oplus h^{2d}(X),\eqno (2)$$
where $h^{i}(X):=h(X,\pi_{i}^{\text{num}},0)\in
\mathscr{M}_{k}^{\text{num}}$\  (See [13] for details).

Furthermore, Murre ([15]) expected that the conjectural
decomposition (1) exists even in $\text{CH}^{d}(X\times
X;\mathbb{Q}):=\text{CH}^{d}(X\times X)\otimes_{\mathbb{Z}}
\mathbb{Q}:=A_{\text{rat}}^{d}(X\times
X)\otimes_{\mathbb{Z}}\mathbb{Q}$ and hence in the category of Chow
motives $\mathscr{M}_{k}^{\text{rat}}$ (w.r.t. the rational
equivalence), $h(X)\in \mathscr{M}_{k}^{\text{rat}}$ has a
decomposition as in (2). In this new setting, the decomposition is
not canonical any more. However, from this conjectural
decomposition, Murre ([15]) conjectured a very interesting
filtration on rational Chow groups which relates the rational
equivalence to the homological equivalence in finite steps as done
by the conjectural Bloch-Beilinson filatration.

More precisely, as in [15], we will say that $X$ has a \emph{Chow-K\"{u}nneth decomposition}
over $k$ if there exist $ \pi_{i}\in \text{CH}^{d}(X\times
X;\mathbb{Q}), \, 0\leq i \leq 2d,$ satisfying

(i) $\pi_{i}$ are mutually orthogonal projectors;

(ii) $\sum_{i}\pi_{i}=\text{cl}_{X\times X}(\Delta_{X});$

(iii) $\text{cl}_{X\times X}(\pi_{i})=\pi_{i}^{\text{hom}}$ (the $i$-th K\"{u}nneth component).

\noindent Equivalently, the Chow motive of $X$ has a
(Chow-K\"{u}nneth) decomposition
$$h(X)=h^{0}(X)\oplus h^{1}(X)\oplus \ldots \oplus h^{2d}(X),$$
where $h^{i}(X):=h(X,\pi_{i},0)\in
\mathscr{M}_{k}^{\text{rat}}.$

Then, Murre proposed in [15] the following famous
conjecture.

\noindent{\bf Murre's Conjecture}

{(A)}: There exists a Chow-K\"{u}nneth decomposition for every irreducible variety
$X\in \mathscr{V}_{k}$ of dimension $d$.

{(B):}
 $\pi_{0},\ldots,\pi_{j-1}$ and $\pi_{2j+1},\ldots,\pi_{2d}$
act as zero on $ \text{CH}^{j}(X;\mathbb{Q}).$

(C):} Let
$F^{v}\text{CH}^{j}(X;\mathbb{Q})=\text{Ker}\pi_{2j}\bigcap
\text{Ker}\pi_{2j-1}\bigcap ...\bigcap \text{Ker}\pi_{2j-v+1}$. Then
the filtration

 \ \ \ \ \ \ $F^{\bullet}$ is independent of
 of the ambiguity in the choice of the $\pi_{i}.$

{(D):}
$F^{1}\text{CH}^{j}(X;\mathbb{Q})=\text{CH}^{j}_{\text{hom}}(X;\mathbb{Q}):=\text{Ker}(\text{cl}_{X}).$

It was shown by Jannsen ([9]) that Murre's conjecture is equivalent
to the Bloch-Beilinson conjecture on rational Chow groups and the
two conjectural filtrations proposed respectively by Murre, and
Bloch and Beilinson  coincides. The main advantage of Murre's
conjecture over Bloch-Beilinson's is that one can check the
statements for specific varieties as we will do in this paper.

Until now, Murre's conjecture is verified for only a few special varieties. It is known that  (A) is true for curves, surfaces ([14]), Abelian
varieties ([17][4]), Brauer-Severi varieties, some threefold
([2][3]), some special fourfold  ([11]), certain modular varieties ([6][7]) and varieties whose Chow motives are finite-dimensional ([10]). As for the other parts of Murre's conjecture, it is known that (B) and (D) are true for the
product of a curve and a surface ([15]), that (B) is true for the
product of two surfaces ([11]) and that some part of  (D) is true
for the product of two surfaces ([11][12]). Jannsen ([10]) proved
that (A), (B), (C) and (D) are true for some very special higher
dimensional varieties over some special ground fields, in particular, he
proved that if $k$ is a rational or elliptic function field (in one
variable) over a finite field $\mathbb{F}$ and $X_{0}$ is an
arbitrary product of rational and elliptic curves over $\mathbb{F}$,
then  (A)-(D) hold for $X_{0}\times_{\mathbb{F}}k$. Gordon and
Murre ([8]) proved that (A)-(D) are true for elliptic modular
threefold over a field of characteristic $0.$

In this paper, we consider Murre's conjecture for certain product
varieties. Concretely, we consider such a problem: if the conjecture
is true for  $X$, when  is it  also true for the product of $X$ with
a curve or some other variety ?  In section 2, we consider the case
of the product of a variety with a projective space. In section 3,
we consider the case of the product of a variety with a curve. In
particular, we generalize Murre's  discussion given in [16], and as
consequences, we prove that if $C$ is a (smooth projective
connected) curve, then Murre's conjecture (B) is true for $X\times
C,$ where $X$ is an elliptic modular threefold over $k$ (an
algebraically closed field of characteristic $0$) or $X$ is an
abelian variety of dimension 3. This implies particularly that (B)
is true for two new kinds of fourfolds other than products of two
surfaces considered in [11] and [12].

\bigskip

\bigskip

{\centerline{\bf 2. Products with projective spaces}}

\bigskip

\noindent Fix a field $k$. Let $X$ (resp. $C$) be a smooth
projective irreducible variety (resp. curve) over $k$. Let $X$ be of
dimension $d$.  In the following, we will always denote by $Z\in
\text{CH}^{j}(X)$ a cycle class. In addition, we denote by $p$ with
some lower indices the projection from a product variety to the
corresponding factors.

In the proof of Theorem 2.3, the following lemma is crucial.

\noindent{\bf Lemma  2.1} ([5]) {\it Let $\mathcal{E}$ be a vector
bundle of rank $r=e+1$ on a scheme $X$ of finite type over}
$\text{Spec}(k)$, {\it with the projection $\pi:
\mathcal{E}\rightarrow X$. Let $\mathbb{P}(\mathcal{E})$ be the
associated projective bundle, $p$ the projection from
$\mathbb{P}(\mathcal{E})$ to $X$, and
$\mathcal{O}_{\mathbb{P}(\mathcal{E})}(1)$ the tautological line bundle
on $\mathbb{P}(\mathcal{E})$. Then there are canonical isomorphisms}
$$\bigoplus_{i=0}^{e}\text{CH}^{j-i}(X) \longrightarrow \text{CH}^{j}(\mathbb{P}(\mathcal{E}))$$
$$(\alpha_{i})\mapsto
\sum_{i=0}^{e}c_{1}(\mathcal{O}_{\mathbb{P}(\mathcal{E})}(1))^{i}\cap
p^{*}\alpha_{i},$$ {\it where}
$c_{1}(\mathcal{O}_{\mathbb{P}(\mathcal{E})}(1))$ {\it is the first
Chern class.}

 $\hfill \Box$

Applying Proposition 3.1 in [5], it is easy to show that the inverse
of the map in Lemma 2.1  is the map
$$\text{CH}^{j}(\mathbb{P}(\mathcal{E})) \longrightarrow \bigoplus_{i=0}^{e}\text{CH}^{j-i}(X),\,\,\beta \mapsto (\beta_{i}),$$
 where $\beta_{e}=p_{*}\beta$  and for $0\leq i \leq e-1$,
$$\beta_{i}=p_{*}(c_{1}(\mathcal{O}_{\mathbb{P}(\mathcal{E})}(1))^{e-i}\cap
\beta-\sum_{t=1}^{e-i}c_{1}(\mathcal{O}_{\mathbb{P}(\mathcal{E})}(1))^{e+t}\cap
p^{*}\beta_{i+t}).$$

\noindent{\bf Lemma  2.2} ([15]) {\it Assume that $Y_{i}\ (i=1,2)$ are
smooth projective irreducible varieties over $k.$ Let $Y=Y_{1} \times
Y_{2}$. If $Y_{i}\ (i=1,2)$ has a Chow-K\"{u}nneth decomposition, then
$Y$ has also a Chow-K\"{u}nneth decomposition.}

$\hfill\Box$

\noindent{\bf Theorem 2.3} {\it Let $X$ be a
smooth projective irreducible variety of dimension $d$ over $k$. If (A),(B) and (D) are
true for $X$, then they are also true for $X\times \mathbb{P}^{r}.$}

\noindent{\it Proof.}
 Let $c=c_{1}(\mathcal{O}_{\mathbb{P}^{r}}(1)),$ that is, the class of any
hyperplane in $\mathbb{P}^{r}.$ For each $0\leq i\leq r$, set
$$\pi_{2i}=c^{r-i} \times c^{i}, \ \ \  \pi_{2i+1}=0,\ \ \ h^{i}(\mathbb{P}^{r})=(\mathbb{P}^{r}, \pi_{i}).$$
Then it is easy to see that there is the Chow-K\"{u}nneth decomposition
$$h(\mathbb{P}^{r})=\bigoplus_{i=0}^{2r}h^{i}(\mathbb{P}^{r})=\bigoplus_{t=0}^{r}h^{2t}(\mathbb{P}^{r}).$$
 Assume that
$X$ has the Chow-K\"{u}nneth decomposition
$$h(X)=\bigoplus_{i=0}^{2d}h^{i}(X), \ \
h^{i}(X)=(X,\pi_{i}^{\prime}).$$
 Then, $X\times \mathbb{P}^{r}$ has
a Chow-K\"{u}nneth decomposition:
$$h(X\times \mathbb{P}^{r})=\bigoplus_{m=0}^{2(d+r)}h^{m}(X\times \mathbb{P}^{r},\pi_{m}),$$
where $ \pi_{m}=\sum_{p+2q=m}\tau_{*}(\pi_{p}^{\prime}\times
c^{r-q}\times c^{q}).$

On the other hand, from Lemma 2.1, we have the isomorphisms:
$$\phi:\text{CH}^{j}(X\times \mathbb{P}^{r};\mathbb{Q})\rightarrow\bigoplus _{i=0}^{r}\text{CH}^{j-i}(X;\mathbb{Q}), \,\, Z \longmapsto (Z_{i}),$$
where $Z_{i}=p_{1*}(([X]\times c^{r-i})\cdot Z)$ with
 $$
\phi^{-1}((Z_{i}))=\sum_{i=0}^{r}([X]\times c^{i})\cdot
p_{1}^{*}Z_{i}=\sum_{i=0}^{r}Z_{i}\times c^{i},$$ and
$$\varphi: \text{CH}^{d+r}(X\times \mathbb{P}^{r} \times X \times\mathbb{P}^{r};\mathbb{Q})\longrightarrow
\bigoplus_{i=0}^{r}\bigoplus_{t=0}^{r}\text{CH}^{d+r-i-t}(X\times
X;\mathbb{Q}), $$
$$\xymatrix@C=0.5cm{
  \alpha \ar@{{|}->}[rrr] && &(\alpha_{it}) }$$
 where $\alpha_{it}=p_{12*}(([X\times X]\times c^{r-t}\times
c^{r-i})\cdot\tau^{*}\alpha).$ In fact, we have
$$\alpha_{it}=p_{13*}(([X]\times c^{r-t}\times [X])\cdot
p_{123*}(([X\times \mathbb{P}^{r}\times X]\times
c^{r-i})\cdot\alpha))$$
$$=p_{13*}(([X]\times c^{r-t}\times [X]\times c^{r-i})\cdot
\alpha)\quad\quad\quad\quad\quad\quad\quad\quad \ \
$$
$$=p_{12*}(([X\times X]\times c^{r-t}\times c^{r-i})\cdot
\tau^{*}\alpha),\quad\quad\quad\quad\quad\quad\quad\quad$$ where
$\tau$ is the isomorphism exchanging the second and the third factor
of the product variety $X\times X \times \mathbb{P}^{r} \times
\mathbb{P}^{r}$. Clearly, we also have
$\varphi^{-1}((\alpha_{it}))=\sum_{i,t}\tau_{*}(\alpha_{it}\times
c^{t} \times c^{i})$.

Now, define the map
$$\Phi: \text{CH}^{d+r}(X\times \mathbb{P}^{r}\times X \times \mathbb{P}^{r};\mathbb{Q}) \times \text{CH}^{j}(X \times \mathbb{P}^{r};\mathbb{Q})\longrightarrow \text{CH}^{j}(X\times \mathbb{P}^{r};\mathbb{Q}),$$
$$\Phi(\alpha,Z):=\alpha(Z):=p_{34*}(\alpha\cdot(Z\times [X\times \mathbb{P}^{r}])).$$
So we have the following diagram
$$
\xymatrix{
  \text{CH}^{d+r}(X\times \mathbb{P}^{r}\times X \times \mathbb{P}^{r};\mathbb{Q}) \times \text{CH}^{j}(X \times \mathbb{P}^{r};\mathbb{Q}) \ar[d]_{\varphi\times\phi} \ar[r]^{\ \ \ \ \ \ \ \ \ \ \ \ \ \ \ \ \ \ \ \ \Phi}
                & \text{CH}^{j}(X\times \mathbb{P}^{r};\mathbb{Q}) \ar[d]^{\phi}  \\
 \bigoplus_{i=0}^{r}\bigoplus_{t=0}^{r}\text{CH}^{d+r-i-t}(X\times X;\mathbb{Q})\times \bigoplus _{i=0}^{r}\text{CH}^{j-i}(X;\mathbb{Q})   \ar@{>}[r]
                & \bigoplus _{i=0}^{r}\text{CH}^{j-i}(X;\mathbb{Q})             }
$$

Here the lower arrow is defined by the other three.

Note that if  $t+i=r, $ then we have
$$\tau_{*}(\pi_{p}^{\prime}\times c^{t} \times c^{q})\cdot
(Z_{i}\times c^{i})=p_{34*}(\tau_{*}(\pi_{p}^{\prime}\times c^{t}
\times c^{q})\cdot (Z_{i}\times c^{i}\times[X\times
\mathbb{P}^{r}]))$$
$$=p_{34*}\tau_{*}((\pi_{p}^{\prime}\times c^{t}
\times c^{q})\cdot (Z_{i}\times[X]\times c^{i}\times
[\mathbb{P}^{r}]))$$
$$=p_{24*}((\pi_{p}^{\prime}\cdot (Z_{i}\times[X]))\times c^{t+i}\times c^{q} ))\ \ \ \ \ \ \ \ \ \ \ \ \ $$
$$=p_{2*}((\pi_{p}^{\prime}\cdot (Z_{i}\times[X]))\times c^{t+i})\times c^{q}\ \ \ \ \ \ \ \ \ \ \ \ \ \ \ $$
 $$ =\pi_{p}^{\prime}(Z_{i})\times c^{q}.\ \ \ \ \ \ \ \ \ \ \ \ \ \ \ \ \ \ \ \ \ \ \ \ \ \ \ \ \ \ \ \  \ \ \ \ \ \ \ \ \ \ $$
 So we conclude that
$$\tau_{*}(\pi_{p}^{\prime}\times c^{t} \times c^{q})\cdot
(Z_{i}\times c^{i})= \left \{ \begin{array}{c}
\pi_{p}^{\prime}(Z_{i})\times c^{q}, \ \ \  \text{if} \  t+i=r;\\
\ \ 0\, ,\ \ \ \  \ \ \ \ \ \ \ \ \ \text{otherwise. }
\end{array}\right.$$
Hence, we can translate the projectors on $X\times \mathbb{P}^{r}$
to those on  $X$ as follows.
$$\phi(\pi_{m}(Z))=\phi\cdot \Phi\cdot (\varphi\times
\phi)^{-1}(\varphi(\pi_{m}),(Z_{i}))=\phi\cdot \Phi(\pi_{m},
\sum_{l=0}^{r}Z_{i}\times c^{i})$$
$$
=\sum_{i=0}^{r}\sum_{p+2q=m}\phi(\tau_{*}(\pi_{p}^{\prime}\times
c^{r-q} \times c^{q})\cdot (Z_{i}\times c^{i}))\ $$
$$
=\sum_{i=0}^{r}\sum_{p+2q=m,i=q}\phi(\pi_{p}^{\prime}(Z_{i})\times
c^{q})\ \ \ \ \ \ \ \ \ \ \ \ \ \ \ \ \ \ \ $$
 $$=(\pi_{m}^{\prime}(Z_{0}), \pi_{m-2}^{\prime}(Z_{i}),
\ldots, \pi_{m-2r}^{\prime}(Z_{r}))\ \ \ \ \ \ \ \ \ \ \ $$

Now, we can prove that the conjectures are true for $X\times \mathbb{P}^{r}.$

For (B), from $0\leq m\leq j-1$ we have $ m-2i\leq (j-i)-1.$ If
$2j+1\leq m\leq 2(d+r)$, then
$$2(j-i)+1\leq m-2i \Longleftrightarrow 2j+1\leq m.$$
So, by the assumptions on $X$, we see that
$$\pi_{m}(Z)=0, \ \ \ \text{for}\ \  0\leq m\leq j-1\ \  \text{or}\ \  2j+1\leq m.$$

For (D),  suppose that $Z\in {\text{CH}}^{j}_{\text{hom}}(X\times
\mathbb{P}^{r};\mathbb{Q})=\text{Ker}({\text{cl}}_{X\times
{\mathbb{P}}^{r}}).$ Then since  $Z_{i}\in
\text{Ker}(\text{cl}_{X})=\text{CH}_{\text{hom}}^{j-i}(X;\mathbb{Q})=\text{Ker}({\pi}
_{2(j-i)}^{\prime})$ by assumption, we have $Z\in \text{Ker}(\pi_{2j}).$ This
completes the proof of Theorem 2.3.

$\hfill\Box$

\noindent{\it Remark} 2.4  (i) We expect that Theorem 2.3 is also
true for non-trivial projective bundles.

(ii) For (C), we can say nothing yet since from the projectors on
$X$ we can get only one but not all set of projectors on
$X\times\mathbb{P}^{r}$.

\noindent{\bf Corollary 2.5} {\it Let $S_{1}, S_{2}$ be smooth projective surfaces over $k.$ Then conjectures (A) and (B) are true for
$S_{1}\times S_{2} \times \mathbb{P}^{r_{1}}\times \ldots \times
\mathbb{P}^{r_{n}}.$}

\noindent{\it Proof:} From Lemma 2.2 and the main theorem of [14], we know that
conjectures (A) and (B) are true for $S_{1}\times S_{2},$ so the
result follows from Theorem 2.3. $\hfill\Box$

\bigskip

\bigskip

\bigskip

{\centerline{\bf 3. Products with curves}}

\bigskip

\noindent Let $C$ be a smooth projective curve over a field $k$ and $e \in C(k).$ It is well-known (see [18] for details)
that $C$ has the Chow-K\"{u}nneth decomposition
$$h(C)=h^{0}(C)\oplus h^{1}(C)\oplus h^{2}(C),$$
where $ h^{i}(C)=(C,\pi_{i}^{\prime\prime})$ with
$$\pi_{0}^{\prime\prime}=[e\times C],\ \
\pi_{2}^{\prime\prime}=[C\times e],\ \
\pi_{1}^{\prime\prime}=\Delta_{C}-\pi_{0}^{\prime\prime}-\pi_{2}^{\prime\prime}.$$

Assume that the irreducible variety $X\in\mathscr{V}_{k}$ has the Chow-K\"{u}nneth decomposition
$$
h(X)= \bigoplus_{i=0}^{2d}h^{i}(X),\ \
h^{i}(X)=(X,\pi_{i}^{\prime}).
$$
Then, the product variety $X\times C$ has the Chow-K\"{u}nneth decomposition
$$h(X\times C)=\bigoplus_{m=0}^{2(d+1)}h^{m}(X\times C,\pi_{m}),$$
where, explicitely,

$\pi_{0}=\pi_{0}^{\prime}\times [e\times
C],\quad\quad\quad\quad\quad\quad\quad\quad\quad\quad\quad\quad\quad\quad\quad\quad\quad\quad\quad$

$\pi_{1}=\pi_{1}^{\prime}\times [e\times
C]+\pi_{0}^{\prime}\times(\Delta_{C}-[e\times C]-[C\times
e]),\quad\quad\quad\quad \ \ \ \ \ \ \ \ \ \ $
$$\pi_{m}=\pi_{m}^{\prime}\times [e\times
C]+\pi_{m-1}^{\prime}\times(\Delta_{C}-[e\times C]-[C\times
e])+\pi_{m-2}^{\prime}\times[C\times e], \  m\geq 2. \ \ $$

Let
 $$\text{CH}_{\text{alg}}^{j}(X;\mathbb{Q}):=\{Z\in
\text{CH}^{j}(X;\mathbb{Q}): Z \sim_{\text{alg}} 0\},$$ where
$Z \sim_{\text{alg}} 0$ means that $Z$ is algebraically equivalent
to $0$.

In the proof of Theorem 3.3, we need the following computations.

\noindent {\bf Lemma 3.1} {\it For any $Z\in \text{CH}^{j}(X\times
C;\mathbb{Q}),$ we have}

(i)  $(\pi_{m}^{\prime}\times[e\times
C])(Z)=\pi_{m}^{\prime}(Z(e))\times [C];$

(ii) $(\pi_{m}^{\prime}\times[C\times
e])(Z)=\pi_{m}^{\prime}(p_{1*}Z)\times [e].$

\noindent{\it Proof:} \ (i) We have
$$(\pi_{m}^{\prime}\times[e\times
C])(Z)=p_{34*}(\tau_{*}(\pi_{m}^{\prime}\times[e\times C])\cdot
(Z\times[X\times C]))$$
$$=p_{34*}(\tau_{*}(\pi_{m}^{\prime}\times[e])\cdot
(Z\times[X])\times [C]$$
$$=p_{3*}(\tau_{*}(\pi_{m}^{\prime}\times[e])\cdot
(Z\times[X])\times [C]\ $$
$$=\pi_{m}^{\prime}(Z(e))\times [C],\ \ \ \ \ \ \ \ \ \ \ \ \ \ \ \ \ \ \ \ \ \ \ \ $$
where $Z(e)=p_{1*}(Z\cdot (X\times[e])).$ Note that in the last
equality, we have used the following computation.
$$p_{3*}(\tau_{*}(\pi_{m}^{\prime}\times[e])\cdot
(Z\times[X])=p_{3*}(\tau_{*}((\pi_{m}^{\prime}\times[C])\cdot[X\times\times[e]\times
X)\cdot (Z\times[X])$$
$$=p_{3*}((p_{13}^{*}\pi_{m}^{\prime}\cdot [X\times e\times
X])\cdot(Z\times [X]))$$
$$=p_{3*}((p_{13}^{*}\pi_{m}^{\prime}\cdot ((Z\cdot[X\times
e])\times [X]))\ \ \ \ \ $$
$$=p_{2*}p_{13*}((p_{13}^{*}\pi_{m}^{\prime}\cdot ((Z\cdot[X\times
e])\times [X]))$$
$$=p_{2*}(\pi_{m}^{\prime}\cdot p_{13*}((Z\cdot[X\times
e])\times [X]))\ \ \ \ \ $$
$$=p_{2*}(\pi_{m}^{\prime}\cdot (p_{13*}(Z\cdot[X\times
e])\times [X]))\ \ \ \ \ $$
$$=\pi_{m}^{\prime}(Z(e)).\ \ \ \ \ \ \ \ \ \ \ \ \ \ \ \ \ \ \ \ \ \ \ \ \ \ \ \ \ \ \ \ \ \ \ \ $$

(ii) Similar to (i), we have
$$(\pi_{m}^{\prime}\times[C\times
e])(Z)=p_{34*}(\tau_{*}(\pi_{m}^{\prime}\times[C\times e])\cdot
(Z\times[X\times C]))$$
$$=p_{34*}(\tau_{*}(\pi_{m}^{\prime}\times[C])\cdot
(Z\times[X]))\times [e])$$
$$=p_{3*}((\tau_{*}(\pi_{m}^{\prime}\times[C])\cdot
(Z\times[X]))\times [e]\ $$
$$=\pi_{m}^{\prime}(p_{1*}Z)\times [e].\ \ \ \ \ \ \ \ \ \ \ \ \ \ \ \ \ \ \ \ \ \ \ \ \ \ \ $$
$\hfill\Box$

\noindent{\bf Lemma  3.2} {\it Let $Z\in \text{CH}^{j}(X\times
C;\mathbb{Q}).$ Then}

 \noindent(i) $(\pi_{m}^{\prime}\times \Delta_{C})(Z)=p_{23*}(p_{13}^{*}Z
\cdot (\pi_{m}^{\prime}\times [C]));$

\noindent(ii) $ (\text{id}_{X}\times f)^{*}((\pi'_{m}\times
\Delta_{C})(Z))=(\pi'_{m})_{K}((\text{id}_{X}\times f)^{*}Z),$ {\it
where $K=k(C)$ is the function field of $C$, $f: \text{Spec}(K)\longrightarrow C$ is the natural morphism and
${(\pi^{\prime}_{m})}_{K}=\pi_{m}^{\prime}\times
\Delta_{\text{Spec}(K)}.$}

\noindent{\it Proof:}\  (i) Let $\delta_{C}: C\longrightarrow C\times
C$ be diagonal morphism. Then, we have
$$(\pi_{m}^{\prime}\times \Delta_{C})(Z)=p_{34*}(\tau_{*}(\pi^{\prime}_{m}\times \Delta_{C})\cdot(Z\times[X\times C]))\ \ \ \ \ \ \ \ \ \ \ \ \ \ \ $$
$$\ \ \ \ \ \ \ \ \ \ \ \ =p_{34*}(\tau_{*}(\text{id}_{X\times X}\times \delta_{C})_{*}(\pi^{\prime}_{m}\times \Delta_{C})\cdot(Z\times[X\times C]))$$
$$=p_{23*}((\pi^{\prime}_{m}\times \Delta_{C})\cdot(\text{id}_{X\times X}\times \delta_{C})^{*}\tau^{*}Z))\ \ \ $$
$$=p_{23*}(p_{13}^{*}Z
\cdot (\pi_{m}^{\prime}\times [C])).\ \ \ \ \ \ \ \ \ \ \ \ \ \ \ \
\ \ \ \  \ \ $$

(ii) From (i) and the following diagram
$$
\xymatrix{
  X_{K}\times_{K}  X_{K} \ar[d]_{\text{id}_{X}\times \text{id}_{X} \times f} \ar[r]^{\ \ \ p_{X_{K}}} & X_{K} \ar[d]^{\text{id}_{X}\times f} \\
  X\times X\times C \ar[r]^{\ \ p_{23}} & X\times C   }
$$
we have
$$(\text{id}_{X}\times f)^{*}((\pi'_{m}\times\Delta_{C})(Z))=(\text{id}_{X}\times f)^{*}p_{23*}(p_{13}^{*}Z\cdot\pi^{\prime}_{m}\times [C])$$
$$ \ \ \ \ \ \ \ \ \  =p_{{X_{K}}*}(\text{id}_{X}\times \text{id}_{X} \times f)^{*}(p_{13}^{*}Z \cdot \pi^{\prime}_{m}\times [C])$$
$$ \ \ \ \ \ \ \ \ \ \ \ \ \ \ \ \ \ \ \ \ \ \ \ \ \ \ \ \ \ \ \ \ \ =p_{{X_{K}}*}((\text{id}_{X}\times \text{id}_{X} \times f)^{*}p_{13}^{*}Z \cdot (\text{id}_{X}\times \text{id}_{X} \times f)^{*}(\pi^{\prime}_{m}\times [C]))$$
$$\ \ \ \ \ \ \ \ \ \ \ \ \ \   =p_{X_{K}*}(((\text{id}_{X} \times f)^{*}Z\times_{K} X_{K})\cdot
\pi^{\prime}_{m}\times \Delta_{\text{spec}K}) $$
$$=(\pi^{\prime}_{m}\times\Delta_{\text{spec}K})((\text{id}_{X} \times
f)^{*}Z)\ \ \ \ \
 $$
$$=(\pi^{\prime}_{m})_{K}((\text{id}_{X} \times
f)^{*}Z).\ \ \ \ \ \ \ \ \ \ \ \ \ \ \    $$

$\hfill \Box$

Our main theorem is the following

\noindent{\bf Theorem 3.3} {\it Let $k$ be an algebraically  closed
field, $X \in \mathscr{V}(k)$ and $C \in \mathscr{V}(k)$ an
irreducible curve with the function field $K=k(C)$. Assume that
$(A)$ and $( B)$ are true for $X$ and $X_{K}$, and that for any $j$,
$\text{CH}_{\text{alg}}^{j}(X_{K};\mathbb{Q})\subseteq
\text{Ker}({(\pi^{\prime}_{2j})}_{K})$. Then $(A)$ and $( B)$ are
also true for $X\times C.$}

\noindent{\it Proof:} The statement about (A) is obvious. We will consider (B) in the following. Let $Z\in \text{CH}^{j}(X\times
C;\mathbb{Q})$.  Easy computations shows that (B) is true if $Z$
is of the form $Z^{\prime}\times [C]$ with $Z^{\prime}\in
\text{CH}^{j}(X;\mathbb{Q}).$ So, we can assume $Z(e)=0,$ since we
have
$$Z=(Z-Z(e)\times[C])+Z(e)\times [C], $$
$$ (Z-Z(e)\times[C])(e)=Z(e)-Z(e)=0.$$

 Assume that
$$0\leq m \leq j-1\ \ \text{ or} \ \  2j+1\leq m
\leq 2(d+1).$$ Then, for $m\geq 1,$ we have
 $$0\leq m-1 \leq(j-1)-1 \ \ \text{or}\ \ 2(j-1)+1\leq m-1,$$
and for $m\geq 2,$ we have
 $$0\leq m-2 \leq(j-1)-1 \ \ \text{or}\ \ 2(j-1)+1\leq m-1.$$

From Lemma 3.1 and the assumptions on $X$, we have (note that ${p_{1}}_{*}Z \in
\text{CH}^{j-1}(X;\mathbb{Q})$)
$$\ \ \ \ \ \ \ \ \ \ (\pi_{m}^{\prime}\times [e\times
C])(Z)=\pi_{m}^{\prime}(Z(e))\times [C]=0,\ \ \ \ \ $$
$$\ \ \ \ \ \ \ \ \ \ \ (\pi_{m-1}^{\prime}\times [e\times
C])(Z)=\pi_{m-1}^{\prime}(Z(e))\times [C]=0,$$
$$\ \ \ \ \ \ \ \ \ \ (\pi_{m-1}^{\prime}\times [C\times
e])(Z)=\pi_{m-1}^{\prime}(p_{1*}Z)\times [e]=0,$$
$$\ \ \ \ \ \ \ \ \ \ (\pi_{m-2}^{\prime}\times [C\times
e])(Z)=\pi_{m-2}^{\prime}(p_{1*}Z)\times [e]=0.$$
So, the problem is
reduced to prove
$$(\pi_{m-1}^{\prime}\times \Delta_{C})(Z)=0,\ \text{if}\ \ 1\leq m \leq j-1 \  \text {or} \  \ 2j+1\leq m \leq 2(d+1).$$

At first, we show that $(\text{id}\times f)^{*}Z$ is algebraically
equivalent to $0$ on $X_{K}.$ In fact, let $\eta$ be the generic
point of $C,$ that is, $K=k(\eta),$ and let $f_{\eta}:
\text{Spec}(K)\longrightarrow C_{K}$ be the $K$-point defined by
$\eta.$ Denote $\eta_{K}=f_{\eta}( \text{Spec}K).$ Then we have
$$Z_{K}(\eta_{K})={p_{X_{K}}}_{*}(Z_{K}\cdot X\times\eta_{K})={p_{X_{K}}}_{*}(Z_{K}\cdot
(\text{id}_{X}\times f_{\eta})_{*}(X\times \text{Spec}K))$$
$$={p_{X_{K}}}_{*}(\text{id}_{X}\times f_{\eta})_{*}((\text{id}_{X}\times f_{\eta})^{*}Z_{K})=(\text{id}_{X}\times f_{\eta})^{*}Z_{K}$$
$$=(\text{id}_{X}\times f_{\eta})^{*}p_{XC}^{XCK*}(Z)=
(\text{id}_{X}\times f)^{*}(Z).\ \ \ \ \ \ \ \ \ \ \ \ \ \ \ $$
Similarly, let $g_{e}: \text{Spec}(K)\longrightarrow C_{K}$ and $g: \text{Spec}k\longrightarrow C$ be the morphisms both defined by $e$. Denote
$e_{K}=g_{e}( \text{Spec}(K)).$ Then we have
$Z_{K}(e_{K})=((\text{id}_{X}\times g)^{*}(Z))_{K}=Z(e)_{K}=0.$

We claim that
$$(\pi^{\prime}_{m-1})_{K}((\text{id}\times
f)^{*}Z)=0,\ \ \text{if}\ \  1\leq m \leq j-1 \ \ \text{or}\ \  2j+1\leq m \leq 2(d+1).$$
In fact,  if $m=2j+1,$ since $(\text{id}\times f)^{*}Z \in
\text{CH}^{j}(X_{K};\mathbb{Q})$ is algebraically equivalent to $0$, from the
assumption $\text{CH}_{\text{alg}}^{j}(X_{K};\mathbb{Q})\subseteq$ $
\text{Ker}((\pi^{\prime}_{2j})_{K})$ we have
$$(\pi'_{m-1})_{K}((\text{id}\times
f)^{*}Z)=(\pi'_{2j})_{K}((\text{id}\times f)^{*}Z)=0;$$ if $1\leq
m\leq j-1$ or $2j+2\leq m\leq 2(d+1)$, we have $1\leq m-1\leq j-2$
or $2j+1\leq m-1\leq 2d+1,$ so from the assumptions on $X$ we get
$(\pi'_{m-1})_{K}((\text{id}\times f)^{*}Z)=0$ since
$(\text{id}\times f)^{*}Z\in \text{CH}^{j}(X_{K};\mathbb{Q}).$

On the other hand, we have the following well known diagram
$$
\xymatrix{
\text{CH}^{j-1}(X\times (C-U);\mathbb{Q})\ar[r] & \text{CH}^{j}(X\times C;\mathbb{Q}) \ar@{->}[r] \ar[d] &  \text{CH}^{j}(X\times U;\mathbb{Q}) \ar[d] \ar[r] & 0 \\
  & \text{CH}^{j}(X_{K};\mathbb{Q}) \ar @{=}[r] & {\underrightarrow{\lim}}_{U\subseteq C} \text{CH}^{j}(X\times U;\mathbb{Q}) &}
$$
where the left vertical map is $z\mapsto (\text{id}\times f)^{*}z$.
So from Lemma 3.2 (ii), we have
$$ (\text{id}_{X}\times f)^{*}((\pi'_{m-1}\times
\Delta_{C})(Z))=(\pi'_{m-1})_{K}((\text{id}\times f)^{*}Z)=0,$$
 hence
$$(\pi_{m-1}^{\prime}\times
\Delta_{C})(Z)=\sum_{i}Z_{i}^{\prime}\times a_{i},\ \text{with}\ \
Z_{i}^{\prime}\in \text{CH}^{j-1}(X;\mathbb{Q})\ \text{and}\ \ a_{i}\in
\text{CH}^{1}(C;\mathbb{Q}).
$$
 In view of $(\pi_{m-1}^{\prime}\times \Delta_{C})^{2}=\pi_{m-1}^{\prime}\times
\Delta_{C},$
 we conclude that
$$(\pi_{m-1}^{\prime}\times \Delta_{C})(Z)=\sum_{i} (\pi_{m-1}^{\prime}\times \Delta_{C})(Z_{i}^{\prime}\times a_{i})\ \ \ \ \ \ \ \ \ \ \ \ \ $$
$$=\sum_{i} p_{23*}[p_{13}^{*}(Z_{i}^{\prime}\times a_{i}) \cdot (\pi_{m-1}^{\prime}\times [C])]$$
$$=\sum_{i} \pi_{m-1}^{\prime}(Z_{i}^{\prime})\times a_{i}=0.\ \ \ \ \ \ \ \ \ \ \ \ \ \ \ \ \ $$
This completes the proof of the theorem.

$\hfill\Box$

Although the theorem above is restricted, we can deduce several
interesting consequences.

\noindent{\bf Corollary 3.4} {\it If $k$ is an algebraically closed
field of characteristic $0$ and $X$ is an elliptic modular threefold
over $k$, then (A) and (B) are true for $X\times C.$}

\noindent{\it Proof:} It was shown in [8] that Murre's conjecture
holds for an elliptic modular threefold over a field of
characteristic 0. Obviously, conjecture (D) for $X_{K}$ implies the
assumption of Theorem 3.3. So, the result is an immediate
consequence of Theorem 3.3.

$\hfill\Box$

\noindent{\bf Corollary 3.5}
 {\it Assume that   algebraic equivalence and  rational equivalence
coincide on $X$. If $(A)$ and $(B)$ are true for $X$ and $X_{K}$,
then $(A)$ and $(B)$ are also true for $X\times C$.}

$\hfill\Box$

\noindent{\it Remark} 3.6\ \ Cellular varieties satisfy the first hypothesis of
the corollary.

By [1] (see also [4] and [15]), for an abelian variety $X$ of dimension $g$ over any field $k$, we have the following decomposition

$$\text{CH}^{j}(X;\mathbb{Q})=\bigoplus_{s=j-g}^{j}\text{CH}_{s}^{j}(X),$$
where $$\text{CH}_{s}^{j}(X):=\{\alpha\in \text{CH}^{j}(X;\mathbb{Q})|n^{*}\alpha=n^{2j-s}\alpha, \forall \  n\in \mathbb{Z}\}.$$

\noindent{\bf Corollary 3.7} {\it Let $X$ be an abelian variety of dimension at most 4 over an algebraically closed field $k$. Assume that for any $j$,
$\text{CH}_{0}^{j}(X_{K})\bigcap
\text{CH}_{\text{alg}}^{j}(X_{K};\mathbb{Q})=0$. Then $(A)$ and
$(B)$ are true for $X\times C$.}

\noindent{\it Proof:} It follows from [1] that conjecture (B) is
true for an abelian variety of dimension at most 4, equivalently, Beauville's vanishing conjecture holds: $\text{CH}_{s}^{j}(X)=0$ if $s<0$.
By assumption and the fact that the algebraic equivalence is adequate, we have
$$\text{CH }_{\text{alg}}^{j}(X_{K};\mathbb{Q})=\bigoplus_{s=1}^{j}(\text{CH}_{s}^{j}(X_{K})\cap \text{CH}_{\text{alg}}^{j}(X_{K};\mathbb{Q})).$$
On the other hand, by Lemma 2.5.1 of [15], we see that for any $j$,
$$\text{Ker}({(\pi^{\prime}_{2j})}_{K})=\bigoplus_{s=1}^{j}\text{CH}_{s}^{j}(X_{K}).$$
Then we can apply Theorem 3.3 to end the proof.

$\hfill\Box$

\noindent{\it Remark} 3.8\ \ The assumption of Corollary 3.7 is a consequence of a conjecture of Beauville: the restricted cycle map $c_{0}: \text{CH}_{0}^{j}(X)\longrightarrow H^{2j}(X)$ is injective for any $j$.

\noindent{\bf Corollary 3.9} {\it Let $X$ be an abelian variety of
dimension 3 over an algebraically closed field $k$. Then} (A) {\it and} (B){\it are  true for $X\times C$.}

\noindent{\it Proof:} By [1], we know that the restricted cycle map
$$c_{0}:\text{CH}_{0}^{j}(X_{K})\longrightarrow H^{2j}(X_{K})$$
is injective for $j=0, 1, g-1,g $ if $X$ is an abelian variety of
dimension $g.$ So for any $j$,
$$\text{CH}_{0}^{j}(X_{K})\bigcap
\text{CH}_{\text{alg}}^{j}(X_{K};\mathbb{Q})\subseteq
\text{CH}_{0}^{j}(X_{K})\bigcap
\text{CH}_{\text{hom}}^{j}(X_{K};\mathbb{Q})=0.$$  Hence the result
follows from Corollary 3.7.

$\hfill\Box$

\bigskip

\bigskip

\bigskip

{\centerline{\bf REFERENCES}}

\bigskip

\noindent 1. A. Beauville.  Sur l'anneau de Chow d'une
vari\'{e}t\'{e} ab\'{e}lienne.  {\it Math.Ann.} {\bf 273}(1986),

647-651

\noindent 2. P. L. del Angel and S. M\"{u}ller-Stach. Motives of
uniruled 3-folds, {\it Compositio Math.}

{\bf 112}(1998), 1-16

\noindent 3. P. L. del Angel,  S. M\"{u}ller-Stach.  On Chow motives
of 3-folds. {\it Trans. Amer. Math.

Soc.} {\bf352}(2000), 1623-1633

\noindent 4. C. Deninger and J. P. Murre.  Motivic decomposition of
abelian schemes and the

Fourier transform. {\it J. Reine Angew. Math.} {\bf422}(1991),
201-219

\noindent 5. W. Fulton.  Intersection theory. {\it Ergeb. Math.
Grenzgeb.} Springer-Verlag, Berlin

1984

\noindent 6. B. Gordon,  M. Hanamura and  J. P. Murre.  Relative
Chow-K\"{u}nneth projectors for

 modular varieties. {\it J. reine
angew. Math.} {\bf 558} (2003), 1-14

\noindent 7.  B. Gordon,   M. Hanamura and   J. P. Murre.  Absolute
Chow-K\"{u}nneth projectors for

modular varieties. {\it J. reine angew. Math.} {\bf580}(2005),
139-155

\noindent 8. B. Gordon and J. P. Murre.  Chow motive of elliptic
modular surfaces and threefolds.

{\it Math. Inst.,} Univ. of Leiden, Report W 96-16 (1996)

\noindent 9.  U. Jannsen.  Motivic sheaves and filtrations on Chow
groups. In {\it Proc. sympos. Pure

Math.} {\bf55}(1994), Part 1 , 245-302

\noindent 10.  U. Jannsen.  On finite-dimensional motive and Murre's
conjecture. Preprint, 2007

\noindent 11.  B. Kahn, J. P. Murre and  C. Pedrini.  On the
transcendental part of the motive of

a surface. Preprint, 2005

\noindent 12.  K. Kimura.  Murre's conjectures for certain product
varieties. {\it J. Math. Kyoto. Univ.}

{\bf 47}(3)(2007), 621-629

\noindent 13. J. P. Murre.  Lectures on motives. In {\it
Transcendental Aspects of Algebraic Cycles,}

Proceedings of Grenoble Summer School(2001), pages 123-170,
Cambridge University

Press, 2004

\noindent 14.  J. P. Murre. On the motive of an algebraic surface.
{\it J. Reine und angew. Math.}

{\bf409}(1990), 190-204

\noindent 15. J. P. Murre.  On a conjectural filtration on the Chow
groups of an algebraic variety,

 I. The general conjectures and some examples. {\it Indag. Math. N. S.} {\bf 4}(2)(1993),

 189-201

\noindent 16. J. P. Murre.  On a conjectural filtration on the Chow
groups of an algebraic variety,

 II. Verification of the conjectures
for threefolds which are the product on a surface

and a curve. {\it Indag. Math. N. S.} {\bf 4}(2)(1993), 189-201

\noindent 17.  A. M. Shermenev. The motive of an abelian variety.
{\it Funct. Anal.} {\bf8}(1974),

47-53

\noindent 18. A.J. Scholl.  Classical motives, In {\it Proceedings
of Symposia in Pure Mathematics,} {\bf

55}(1994), Part I, 163-187

\bigskip

\bigskip

\bigskip

\bigskip

\noindent {\large K}EJIAN {\large X}U \ \ \ \ \ kejianxu@amss.ac.cn

\bigskip

\noindent College of Mathematics

\noindent Qingdao University

\noindent  Qingdao 266071, China

\bigskip

\bigskip

\noindent {\large Z}E {\large X}U \ \ \ \ \ xuze@amss.ac.cn

\bigskip

\noindent Institute of Mathematics

 \noindent Academy of Mathematics and System Science

 \noindent Chinese Academy of Sciences

 \noindent Beijing 100190, China

\end{document}